\documentclass[12pt]{article}

\newenvironment{pf}{\noindent\textbf{Proof.}\quad}{\hfill{$\Box$}}
\usepackage[colorlinks=false]{hyperref}
\usepackage[dvips]{graphicx}
\usepackage[usenames,dvipsnames]{color}
\usepackage{epsfig}
\usepackage{rotating}
\usepackage{amssymb}
\usepackage{amsmath,amsfonts}
\usepackage[latin1]{inputenc}
\usepackage{verbatim}
\usepackage[tableposition=top,font=small,labelfont=bf,format=hang]{caption}
\usepackage{booktabs}

\newtheorem{lem}{Lemma}
\newtheorem{thm}{Theorem}
\newtheorem{prop}{Proposition}
\newtheorem{conj}{Conjecture}
\newtheorem{rem}{Remark}
\newtheorem{coro}{Corollary}
\newtheorem{ex}{Example}

\newcommand{\F}{\mathbb{F}}
\newcommand{\FF}{\mathbb{F}}

\def\GL{\mathrm{GL}}
\def\Aut{\mathrm{Aut}}
\def\SAut{\mathrm{SAut}}

\DeclareMathOperator{\wt}{wt}
\begin{document}

\title{Self-dual Hadamard bent sequences\thanks{This work is supported in part by the National Natural Science
Foundation of China (12071001), the Excellent Youth Foundation of
Natural Science Foundation of Anhui Province (1808085J20). The work of Dean Crnkovi\'c
is supported by Croatian Science Foundation under the project 6732.}}

%


\author{
Minjia Shi,  Yaya Li\thanks{Minjia Shi and Yaya Li, School of Mathematical Sciences, Anhui University, Hefei, Anhui, 230601,
China, {\tt  smjwcl.good@163.com, yayali187125@163.com }}, \
Wei Cheng\thanks{T\'el\'ecom Paris; Secure-IC S.A.S., 104 Boulevard du Montparnasse, 75014 Paris, France, {\tt wei.cheng@telecom-paris.fr}}, \
Dean Crnkovi\'c\thanks{Faculty of Mathematics, University of Rijeka, Croatia, {\tt deanc@math.uniri.hr}}, \\
Denis Krotov\thanks{Sobolev Institute of Mathematics, Novosibirsk 630090, Russia, {\tt krotov@math.nsc.ru}}, \
Patrick Sol\'e\thanks{CNRS,  University of Aix Marseille, Centrale Marseille, I2M, Marseille, France, {\tt sole@enst.fr}}
}

\date{}

\maketitle
\begin{abstract}
A new notion of bent sequence related to Hadamard matrices was introduced recently, motivated by a security application (Sol\'e et al, 2021). We study the self-dual class in length at most $196$.
We use three competing methods of generation: Exhaustion, Linear Algebra and Groebner bases.
Regular Hadamard matrices and Bush-type Hadamard matrices provide many examples. We conjecture that if $v$ is an even perfect square, a self-dual bent sequence of length $v$ always exists.
We introduce  the strong automorphism group of Hadamard matrices, which acts on their associated self-dual bent sequences. We give an efficient algorithm to compute that group.
\end{abstract}


{\bf AMS Classification (MSC 2010):} Primary 94D10, Secondary 15B34
\section{Introduction}
Bent functions and bent sequences are classical objects in algebraic combinatorics with sundry connections to design theory, distance regular graphs, and symmetric cryptography~\cite{Mesnager:2016:Bent}, \cite{Tokareva:2015}.
In~\cite{SCGR:2021isit} a new notion of bent sequence was introduced as a solution in $X,Y$ to the system
$$\mathcal{H}X=Y, $$
where $H$ is a  Hadamard matrix of order $v$, normalized to $\mathcal{H}=H/\sqrt{v}$ and  $X,Y \in \{\pm 1\}^v$.
Given $H$ the vector $X$ defines a {\it Hadamard bent} (binary) sequence by the correspondence
$$x\mapsto (X_x+1)/2. $$
When $v$ is a power of $2$ and  $H$ is the Hadamard matrix of  Sylvester type,  we recover the classical notion of bent sequence~\cite{Mesnager:2016:Bent}, \cite{Tokareva:2015}. They meet
the covering radius of Hadamard codes (see \S \ref{sec:back}, Lemma \ref{CR}) in the same
way that classical bent functions meet the covering radius of the first order Reed-Muller code \cite[Chap.14,\,Th.6]{MWS}. Beyond generalization for the sake of generalization, this notion was introduced in~\cite{SCGR:2021isit} from a cryptographic perspective (See \S \ref{sec:application} for details).
We believe this concept has a combinatorial interest of its own, as it pertains to the fine print of Hadamard matrix theory: regular matrices, and automorphism groups.

It is proved in~\cite{SCGR:2021isit} that this
kind of bent sequence can only exist if $v $ is a perfect square. As is well-known, Hadamard matrices of order $>2$ only exist for $v$ a multiple of $4$. Thus we reduce to $v=4m^2=(2m)^2$ with $m$ an integer (in practice $1 \le m \le 7$).

In~\cite{CDPS:2010}, when $H$ is of Sylvester type, a linear algebra technique is used to find {\bf self-dual} bent sequences, a situation
which corresponds to the case $X=Y$ in the above equation. Namely $X$ is, in particular, an eigenvector
associated to the eigenvalue $1$ of $\mathcal{H}$. The condition $X \in \{\pm 1\}^v$ has to be checked independently, using a basis of the eigenspace.

In order for the approach of~\cite{CDPS:2010} to work for more general Hadamard matrices than Sylvester type, we need to assume that $\mathcal{H}$ has the eigenvalue  $1$ in its spectrum, and that the dimension of the associated eigenspace
is not too large, as this parameter controls the complexity of the search. Another algebraic technique consists in reducing the existence
of a sequence of length $v$ to a quadratic system in $v$ variables, which can be solved by using Groebner bases. This second method works well as long as the number of variables is less than one hundred.
The brute force approach which consists in checking all $2^v$ possible $X$'s is not feasible for $v>30$, say.

A connection with regular Hadamard matrices is pointed out. Every regular Hadamard matrix admits the all-one vector as a self-dual bent sequence.
In particular,  Bush-type Hadamard matrices of order $v=4u^2$ afford at least $2^{2u}$ self-dual bent sequences. The Karaghani conjecture \cite{Kharaghani:2000:twin} on the existence of regular Hadamard matrices
suggests then that self-dual bent
sequences exist for all even perfect square orders (Conjecture \ref{conjecture}). This conjecture is satisfied for all even square orders where a regular Hadamard matrix exists. The first unknown order seems
to be $v=4u^2$ for $u=47.$

We introduce the notion of strong automorphism group of a Hadamard matrix. This group acts on the associated self-dual bent sequences. We give an efficient
algorithm to compute it, based on a digraph defined from the matrix. We also connect this group
to the group of polarities of the Menon design defined by the matrix.

The material is organized as follows. The next section collects notions and notations needed in the following sections. Section~\ref{sec:prop} investigates interesting properties of self-dual bent sequences, and Section~\ref{sec:search} develops the three search methods. Then Section~\ref{sec:num}
displays the numerical results we found. Section~\ref{sec:application} introduces the application of Hadamard bent sequences. Section~\ref{sec:con} concludes the article. An appendix develops the construction techniques of Hadamard matrices of order $16,\,36,\,64,\,100,\,144$ and $196$.
\section{Background material}
\label{sec:back}

\subsection{Hadamard matrices}
A {\bf Hadamard matrix} $H$ of order $v$ is a $v$ by $v$ real matrix with entries $\pm 1$ satisfying $HH^{t}=vI_{v}$, where $H^{t}$ is the transpose of $H$ and $I_{v}$ is the identity matrix of order $v$. If $v>2$ it is well-known that $v$ must be a multiple of $4$~\cite{MWS}.
An important construction of Hadamard matrices, due to {\bf Sylvester} is obtained for $v=2^h$, when $H$ is indexed by binary vectors of length $h$ and $H_{xy}=(-1)^{\langle x,y\rangle}$,
where ${\langle x,y\rangle}=\sum\limits_{i=1}^hx_iy_i$. We denote henceforth  this matrix by $S_v$. For more general constructions, properties and applications of Hadamard matrices we refer the reader to~\cite{Horadam:2007}.
A Hadamard matrix of order $v$ is {\bf normalized} if both first row and first column are equal to the all-one vector.
A Hadamard matrix of order $v$ is {\bf regular} if its $v$ row and column sums are all equal to a constant $\sigma$. In that case, it is known that $v=4u^2$ with $u$ a positive integer and that
$\sigma=2u$ or $-2u$~\cite{WSW:72}. A special class is that of {\bf Bush-type} Hadamard matrices~\cite{Bush:71}. A Hadamard matrix of order $v=4u^2$ is said to be { Bush-type}
if it is blocked into $2u$ blocks of side $2u,$
denoted by $H_{ij},$ such that the diagonal blocks $H_{ii}$ are all-ones, and that the off-diagonal blocks have row and column sums zero.
\subsection{Bent Boolean functions}
A Boolean function $f$ of arity $h$ is any map from $\F_2^h$ to $\F_2$.
The sign function of $f$ is defined by $F(x)=(-1)^{f(x)}$.
The Walsh-Hadamard transform of $f$ is defined as
$$ \widehat{f}(y)=\sum_{x \in \F_2^h} (-1)^{{\langle x,y\rangle}+f(x)}. $$
Thus in term of vectors $$ \widehat{f}=S_v F. $$
A Boolean function $f$ is said to be {\bf bent} iff its Walsh-Hadamard transform takes its values in $\{ \pm 2^{h/2}\}$. Such functions can only exist if $h$ is even.
The dual of a bent function $f$ is defined by its sign function $ \widehat{f}/2^{h/2}$~\cite{CDPS:2010}.
A bent function is said to be {\bf self-dual} if it equals its dual.
In terms of the Sylvester matrix the sign function $F$ of a self-dual bent function satisfies $\mathcal{S}_v F=F$ where $\mathcal{S}_v=\frac{S_v}{2^{h/2}}$ and $v=2^h$.
\subsection{Hadamard bent sequences}

If $H$ is a Hadamard matrix of order $v$ a {\bf bent sequence} of length $v$ attached to $H$ is any vector $X \in \{\pm 1\}^v$, such that
$$\mathcal{H}X=Y, $$
where  $\mathcal{H}=H/\sqrt{v}$ and  $Y \in \{\pm 1\}^v$.

The {\bf dual} sequence of $X$ is defined by $Y=\mathcal{H}X$. If $Y=X$, then $X$ is a {\bf self-dual} bent sequence attached to $H$. It is easy to see that the vector $Y$ is itself a bent sequence attached to $H^t$.

When $H=S_v$ we recover the definitions of the preceding subsection.
\subsection{Hadamard codes}

We consider codes over the alphabet $A=\{ \pm 1\}$. If $H$ is a Hadamard matrix of order $v$, we construct a code $C$ of length $v$ and size $2v$ by taking
the columns of $H$ and their opposites. Let $d(.,.)$ denote the Hamming distance on $A$.
The {\bf covering radius} of a code $C$ of length $v$ over $A$ is defined by the formula
$$r(C)= \max_{y \in A^v}\min_{x \in C} d(x,y). $$

The following lemma is immediate by Theorems 1 and 2 of~\cite{SCGR:2021isit}.
{\lem \label{CR} Let $v$ be an even perfect square, and let $H$ be a Hadamard matrix of order $v$, with the associated Hadamard code $C$.
The vector $X \in A^v$ is a bent sequence attached to $H$ iff
$$\min_{Y \in C} d(X,Y)=r(C)=\frac{v-\sqrt{v}}{2} . $$}

This Lemma generalizes nicely Theorem 6 of \cite[Chap. 14]{MWS}.
\subsection{Graphs}
A directed graph ({\bf digraph} for short)  on a set $V$ of vertices is determined by a set of arcs $E \subseteq V\times V$.
Declare two vertices $x,y$ adjacent and write $ x\sim y$  iff $(x,y) \in E$. The {\bf adjacency matrix} $A$ is then  defined by
$$A_{xy} = \begin{cases}
           1 \,\quad \mbox{if} \,   x\sim y, \\
           0 \,\quad \mbox{if} \,  x\nsim y.
          \end{cases}$$

The {\bf automorphism group} of such a digraph is the group of permutations on $V$ that preserve incidence.

\section{ Properties of self-dual bent sequences}
\label{sec:prop}
\subsection{Automorphism groups}
The class of Hadamard matrix of order $v$ is preserved by the three following operations:
\begin{itemize}
 \item row permutation,
 \item column permutation,
 \item row or column negation,
\end{itemize}
which form a group $G(v)$ with structure $(S_v \wr S_2)^2$, where $S_m$ denotes the symmetric group on $m$ letters. We denote by $S(v)$ the group of diagonal matrices of order $v$  with diagonal elements in $\{\pm 1\}$,
and by $M(v)$ the matrix group generated by $P(v)$, the group of {\bf permutation matrices} of order $v$,  and $S(v)$.
The action of $G(v)$ on a Hadamard matrix $H$ is of the form
$$ H\mapsto PHQ, $$ with $P,Q \in M(v)$.
The {\bf automorphism group} $\Aut(H)$ of a Hadamard matrix $H$ is defined classically as the set of all pairs $(P,Q) \in G(v)$ such that $PHQ=H$~\cite{Hall:62:Mathieu}.
Some information on this group in the case of Paley matrices can be found
in~\cite{CES:2018},
\cite{Kantor:69:aut},
\cite{deLauFla:2011}.
The cases of type I and type II
(i.e., $q\equiv 3 \pmod{4}$ and $q\equiv 1 \pmod{4}$ where $q$ is the prime power in the definition of Paley matrix) were exactly determined in~\cite{Kantor:69:aut} and~\cite{deLauSta:2008}, respectively.
The automorphism group of (a generalization of) the Sylvester matrix can be found in~\cite[p.\,101--103]{deLauFla:2011}.
We give a characterization of $\Aut(H)$
for the Sylvester matrix $S_v$.

Consider the action of an {\bf extended
affine transform} $T_{A,b,d,c}$ on a Boolean function $f$, i.e.,
$$f(x) \mapsto f(A^{-1}x+A^{-1}b)
 \cdot (-1)^{{\langle d,x\rangle}}
 \cdot c, $$
where $A$ is an $m$-by-$m$ invertible matrix over $\FF_2$,
$b \in \FF_2^m$,
$d \in \FF_2^m$,
$c\in\{1,-1\}$.

{\thm The pair $(T_{A,b,d,c},T_{(A^{-1})^t,d,b,c(-1)^{\langle b,d\rangle}})$
 is in $\Aut(S_v)$. }

\begin{pf}
 By definition of $S_v$, with $v=2^m$, we have $g=S_v f$ iff $g(y) =
\sum\limits_{x \in \FF_2^m}f \cdot (-1)^{\langle x,y\rangle}.$ We compute $S_v T_{A,b,d,c} (f)$ by the same formula.

\begin{equation*}
\begin{split}
	\sum_{x \in \FF_2^m} T_{A,b,d,c} &f(x)\cdot (-1)^{\langle x,y\rangle} \\
	&= \sum_{x \in \FF_2^m}
	f(A^{-1}x+ A^{-1}b) \cdot (-1)^{\langle d,x\rangle} \cdot c \cdot (-1)^{\langle x,y\rangle} \\
	&= c\cdot \sum_{x' \in \FF_2^m}
	f(x') \cdot (-1)^{\langle Ax'+b,y+d\rangle} \quad \textit{\small//By taking $x=Ax'+b$}\\
	&= c\cdot (-1)^{\langle b,y+d\rangle}\cdot \sum_{x' \in \FF_2^m}
	f(x') \cdot (-1)^{\langle Ax',y+d\rangle}\\
	&= c(-1)^{\langle b,d\rangle}\cdot (-1)^{\langle b,y\rangle}\cdot \sum_{x' \in \FF_2^m}
	f(x') \cdot (-1)^{\langle x',A^ty+A^td\rangle}\\
	&= T_{(A^{-1})^t,d,b,c(-1)^{\langle b,d\rangle}}
	g(y)
\end{split}
\end{equation*}
Thus $S_vT_{A,b,d,c}f=T_{(A^{-1})^t,d,b,c(-1)^{\langle b,d\rangle}}
 g$, and the pair $(T_{A,b,d,c}, T_{(A^{-1})^t,d,b,c(-1)^{\langle b,d\rangle}})$
 is in $\Aut(S_v)$.
\end{pf}

To work on the symmetries
of bent sequences we will require the notion of
{\bf strong  automorphism group} $\SAut(H)$ of $H$ defined as the set of $P \in M(v)$ such that $PH=HP$. Then we can state the following result.

{\prop If $X$ is self-dual bent sequence
for $H$, and if $P\in M(v)$ is a strong automorphism of $H$, then $PX$ is also self-dual bent sequence for $H$.}

\begin{pf}
 By hypothesis $\mathcal{H}X=X$. Multiplying on left this equation by $P$ we get
 $$PX=P \mathcal{H}X=\mathcal{H}PX. $$ Letting $Y=PX$, we see that $\mathcal{H}Y=Y$.
\end{pf}

A partial characterization in the case of $\SAut(S_v)$ is as follows. It is an immediate corollary of the preceding theorem and its proof is omitted.

{\coro An extended affine transform $T_{A,b,d,c}$ is in $\SAut(S_v)$ iff $A^t = A^{-1}$, $b=d$ and $\wt_H(b)$ is even.}
\begin{rem}
{\em In particular, the number of such transforms is $|\mathcal{O}_m|2^m$ where $\mathcal{O}_m=\{ A \in \GL(m,\F_2) \mid AA^t=I\}$. By~\cite[Theorem 4]{Janusz:2007}, we know that}
\end{rem}
\begin{itemize}
 \item $|\mathcal{O}_m|=2^{k^2}\prod\limits_{i=1}^{k-1}(2^{2i}-1)$ if $m=2k$,
 \item $|\mathcal{O}_m|=2^{k^2}\prod\limits_{i=1}^{k}(2^{2i}-1)$ if $m=2k+1$.
\end{itemize}
For the first few values of $m$, we get $1,
2,
8,
48,
768,
23040,
1474560,
185794560$.

A stronger characterization of the automorphism group of the set of self-dual bent functions within all Hamming isometric maps is in~\cite{Kutsenko:2020}.
A weaker form of our corollary appears in~\cite[Theorem 1]{FSSW:2013} where the group of invertible matrices $A$ satisfying $A^{t}=A^{-1}$ is called the orthogonal group.
An algorithm to compute the strong automorphism group is given at the end of the section.

Two Hadamard matrices $H$ and $K$ are {\bf strongly equivalent} if there is $P\in M(v)$  such that $PHP^t=K$. (This relation is an equivalence relation on the set of Hadamard matrices). Then they share the same self dual bent sequences, up to a monomial transform, as the next result, the main motivation for this new concept, shows.

{\prop If $H$ and $K$ are {strongly equivalent} Hadamard matrices, satisfying $K=PHP^t$, with $P\in M(v)$ then their respective sets of self-dual bent sequences, say $S(H)$
and $S(K)$, satisfy $S(H)=P^{t}S(K)$. }

\begin{pf}
 If $KX=\sqrt{v}X$ for some $X\in \{\pm 1\}^v$, then let $Y=P^tX$. We see that $HY=\sqrt{v}Y$ and that $Y\in \{\pm 1\}^v$. The result follows.
\end{pf}

The database of Magma collects the orbits of Hadamard matrices under $G(v)$ by their normalized representative. It is plain to see that the action of $G(v)$ does not preserve
the self-dual bentness property. A simple example is given by the pair of equivalent matrices $H$ and $-H$ who cannot allow a common nonzero self-dual bent sequence. In fact, the action of $G(v)$
can produce self-dual bent sequences as the next result shows.

{\prop If $X$ is bent sequence for $H$, then there is an equivalent Hadamard matrix $H'$ such that $X$ is self-dual bent sequence for $H'$.}

\begin{pf}
 By hypothesis $\mathcal{H}X=Y$. There is a  matrix $S \in S(v)$  such that $Y=SX$. Since $S$ is an involution we have
 $$ X=S\mathcal{H}X=\mathcal{H'}X, $$ where $H'=SH$.
\end{pf}

\subsection{Regular Hadamard matrices}
A direct connection between Hadamard bent sequences and regular Hadamard matrices is as follows.

{\prop If $H$ is a regular Hadamard matrix of order $v=4u^2$, with $\sigma=2u$, then $j$ is a self-dual bent sequence for $H$ where j is the all-one vector of length $v$.}

\begin{pf}
 By definition of regular Hadamard matrices $Hj=\sigma j=\sqrt{v}j$, yielding $\mathcal{H}j=j$.
\end{pf}

Any construction of regular Hadamard matrices implies the existence of self-dual Hadamard bent sequences. The reference~\cite{Crn:2006} yields the following result.

{\coro Let $p$ and $2p-1$ be prime powers and $p\equiv 3 \pmod{4}$, then there exists a self-dual Hadamard bent sequence of length $4p^2$.
In particular $p=3$ yields a self-dual Hadamard bent sequence of length $36$, and $p=7$ yields a self-dual Hadamard bent sequence of length $196$.}

Another construction, valid for some primes $\equiv 7 \pmod{16}$ can be found in~\cite{LMS:2006}.

In fact each Bush-type Hadamard matrix implies the existence of many self-dual bent sequences.

{\prop If $H$ is a Bush-type  Hadamard matrix of order $v=4u^2$, then there are at least $2^{2u}$ self-dual bent sequences for $H$.}

\begin{pf}
 From the definition, we see that the sequence $X$ defined by  $$X^t=(\pm j,\dots,\pm j), $$  where $j$ is the all-one vector of length ${2u}$, and the $2u$ signs $\pm1$ are arbitrary, is self-dual bent sequence.
\end{pf}

In view of Kharagani's conjecture that Bush-type Hadamard matrices exist for all even perfect square orders~\cite{Kharaghani:2000:twin}, the two previous propositions suggest the following.

{\conj \label{conjecture} If $v$ is an even perfect square, then there exists a self-dual Hadamard bent sequence for some Hadamard matrix of order $v$.}

We show that the Kronecker product of two self-dual bent sequences is also self-dual bent sequence. Recall that the {\bf Kronecker product} $K=X \otimes Y$ of two sequences $X$ and $Y$
of respective lengths $v$ and $w$ is defined by $K_{(i,j)}=X_iY_j$. Similarly, the {\bf Kronecker product} of two Hadamard matrices $U$ and $V$ of respective orders $v$ and $w$ can be defined as
$$(U \otimes V)_{(i,j),(k,\ell)}=U_{ij}V_{k\ell}. $$

{\prop If $X$ and  $Y$ are two self-dual bent sequences with respective Hadamard matrices $U$ and $V$, then $X \otimes Y$ is a self-dual bent sequence attached to $(U \otimes V)$.}

\begin{pf}
 As is well-known~\cite{WSW:72}, if both $U$ and $V$ are Hadamard matrices then so is $(U \otimes V)$. Now the relations $X \sqrt{v}=U X$ and $Y \sqrt{w}=V Y$ entail
 $$(U \otimes V)(X \otimes Y)=\sqrt{v w}(X \otimes Y). $$
 This completes the proof.
\end{pf}

This implies for instance, the existence of self-dual bent sequences of length $64=4 \times 16$, from the existence of self-dual bent sequences in lengths $4$ and $16$.

\subsection{Computing the strong automorphism group}
\subsubsection{The strong automorphism group}

Define a digraph $G(H)$ by its adjacency matrix $A(H)$ as follows.
This matrix is obtained  by replacing in $H$
\begin{itemize}
\item
the  $1$'s by
$\left[\begin{array}{cc}1&0\\0&1\end{array}\right]$,
\item
the  $-1$'s by
$\left[\begin{array}{cc}0&1\\1&0\end{array}\right]$.
\end{itemize}

{\thm The group $\SAut(H)$ is isomorphic to the automorphism group of $G(H)$. }

\begin{pf}
 First, we note that any automorphism
of $G(H)$ do not break the \emph{blocks}
$\{0,1\}$, $\{1,2\}$, \ldots, $\{2n-2,2n-1\}$
(we assume that the vertices of $G(H)$
are the indices of the corresponding columns/rows in $A(H)$).
Indeed, two vertices are in the same block
if and only if their neighborhoods do not intersect.

 The rest is straightforward.
Permuting blocks of vertices in $G(H)$ corresponds to
permuting the column/row indices of $H$,
while swapping two vertices in the same block
corresponds to the negation of the corresponding
row and column in $H$.
\end{pf}
\begin{rem}
{\em This graphical method can also be used to check if two Hadamard matrices are strongly equivalent.}
\end{rem}
\subsubsection{The permutation part}
The permutation part $C(H)$ of the strong group defined by $$C(H)=\{ P \in P(v) \mid PH=HP \}$$ admits an intuitive interpretation in terms of directed graphs ({\bf digraphs}).
Let $\Gamma(H)$ denote the digraph with adjacency matrix $A$ where $H=J-2A$, and $J$ denote the $v$ by $v$ all-one matrix.

{\thm The group $C(H)$ is the group of isomorphisms of  $\Gamma(H)$.}

\begin{pf}
 Since $P\in P(v)$, we have $PJ=JP=J$.
 Thus $HP=PH$ iff $PA=AP$. Assume now that  $i,j$ have respective preimages  $h$ and $k$ under $P$. Computing matrix products we get
$$(PA)_{hj}=a_{ij}=(AP)_{hj}=a_{hk}. $$
Thus $ i  \sim j$ iff $ h  \sim k$ which shows that $P$ preserves adjacency in $\Gamma(H)$ .
\end{pf}
\begin{rem}
{\em The above proof is a direct extension of the proof of~\cite[Prop.\,15.2]{Biggs:74} from graphs to digraphs.}
\end{rem}
\begin{ex}
{\em Let $H$ be the Paley type II Hadamard matrix of order $36$. Then Magma commands {\tt HadamardAutomorphismGroup} and {\tt AutomorphismGroup} allow us to compute}
\end{ex}
\begin{itemize}
 \item $|\Aut(H)|=2^7\times 3^2 \times 17$,
 \item $|\SAut(H)| = 2^5 \times 3^2 \times 17$,
 \item $|C(H)|=2^3 \times 17$.
\end{itemize}

Note that the latter number divides the former, as it should, since $C(H)$ can be embedded in a subgroup of $\Aut(H)$ by writing $PHP^t=H$.
More generally $ |C(H)|$ divides $|\SAut(H)| $ which divides $|\Aut(H)|$.
In the next table we give the same information for the $5$ matrices of order $16$ in Magma database. The first row is the index $j$ of $H$ in the Magma database.
$$\begin{array}{|c|c|c|c|c|c|} \hline
j&1&2&3&4&5\\ \hline
|\Aut(H)|& 2^{15} \times 3^2 \times 5 \times 7& 2^{12} \times 3 \times 7&2^{12} \times 3 \times 7& 2^{15} \times 3^2&2^{14} \times 3\\ \hline
|\mathrm{SAut}(H)|  & 2^9 \times 3^2 \times 5 & 2^2 & 2 & 2 & 2 \\ \hline

 |C(H)|& 2^4 \times 3^2 \times 5& 1&1&1&1\\ \hline

\end{array}$$

\subsubsection{Involutions}

Define further the group $C_2(H)=\{ P \in C(H) \mid P^2=I\}$, consisting of the identity and of the involutions in $C(H)$.
This can be interpreted in terms of {\bf combinatorial designs}.
Consider the incidence system $(\mathcal{V}, \mathcal{B}, \mathcal{I})$
defined by the following three rules:
\begin{itemize}
 \item $\mathcal{V}$ is the set of rows of
$H$,
 \item $\mathcal{B}$ is the set of columns of $H$,
 \item $i\, \mathcal{I} \, j$ iff $H_{ij}=-1$.
\end{itemize}

A {\bf duality} $\pi$ of this incidence system on its dual $(\mathcal{B}, \mathcal{V},  \mathcal{I})$ is then defined as a bijection $\pi$ between
$\mathcal{V}$ and $\mathcal{B}$ that
preserves incidence: $\mathcal{B}$ and $\mathcal{V}$ are swapped by $\pi$ and  $i \mathcal{I} j$ iff $\pi(i) \mathcal{I} \pi(j) $ (Cf~\cite[(4.1.b) p.\,34, Def.\,4.9]{BJH:99}).
The set of all dualities form a group for map composition. Furthermore if a duality is an involution, it is called  a  {\bf polarity}. In terms of the incidence matrix $A$ of $\mathcal{I}$,
a permutation matrix $P$ is a polarity if $PA=A^tP^t$ and $P^t=P$, or, equivalently, if $PA=A^tP$ and $P^t=P$.

{\thm If $H$ is symmetric, then the  group $C_2(H)$ coincides with the  group of polarities of the above incidence structure.}

\begin{pf} The incidence matrix $A$ of $(\mathcal{V}, \mathcal{B}, \mathcal{I})$ satisfies by definition $H=J-2A$, where $J$ denotes the $v$ by $v$ all-one matrix. If $H$ is symmetric, then $H=H^t$, and, since $J=J^t$, we have $A=A^t$.
 Since $P\in P(v)$, we have $PJ=JP=J$. Thus $HP=PH$ iff $PA=AP$, or, equivalently, iff $PA=A^tP$. The result follows.
\end{pf}

\section{Search Methods}
\label{sec:search}

\subsection{Exhaustion}

This method is only applicable for small $v$'s.

\begin{enumerate}
 \item[\rm(1)] Construct $H$ a  Hadamard matrix of order $v$ like in~\cite{SCGR:2021isit} by using Magma database.
 \item[\rm(2)] For all $X \in \{\pm 1\}^v$ compute $Y=\mathcal{H}X$. If $Y=X$, then $X$ is self-dual bent sequence attached to $H$.
\end{enumerate}
{\bf Complexity:} Exponential in $v$ since $| \{\pm 1\}^v|=2^v$.

\subsection{Linear algebra}
\label{sec:search:lin}

This method is more complex to program than the others but allow to reach higher $v'$s.
\begin{enumerate}

\item [\rm(1)] Construct $H$ a Hadamard matrix of order $v.$ 
Compute $\mathcal{H}=\frac{1}{\sqrt{v}}H.$
\item [\rm(2)] Compute a basis of the eigenspace associated to the eigenvalue $1$ of $\mathcal{H}$.
\item [\rm(3)] Let $B$ denote a matrix with rows such a basis of size $k \le  v$. Pick $B_k$ a $k$-by-$k$  submatrix of $B$ that is invertible, by the algorithm given below.
\item [\rm(4)] For all  $Z \in \{\pm 1\}^k$ solve the system in $C$ given by $Z=CB_k$.
\item [\rm(5)] Compute the remaining $v-k$ entries
 of $CB$.
\item [\rm(6)] If these entries are in $\{\pm 1\}$ declare $CB$  a self-dual bent sequence attached to $H$.
\end{enumerate}

To construct $B_k$ we apply a greedy algorithm. We construct the list $J$ of the indices of the columns of $B_k$ as follows.

\begin{enumerate}
 \item [\rm(i)] Initialize $J$ at $J=[1]$.
 \item [\rm(ii)] Given a column of index $\ell$ we compute the ranks $r$ and $r'$ of the submatrices of $B$ with $k$ rows and columns defined by the respective lists $J$ and $J'=Append(J,\ell)$.
 \item [\rm(iii)] If $r<r'$ then update $J:=J'$.
 \item [\rm(iv)] Repeat until $|J|=rank(B)$.
\end{enumerate}
\begin{rem}
{\em 

  The matrix in step 1 can be constructed by using Magma database or by the techniques in the Appendix.
  }
\end{rem}

  \begin{rem}
  {\em 
 If the first column of $B$ is zero, step $\rm(i)$ does not make sense, but then there is no self-dual bent sequence in that situation, as all eigenvectors have first coordinate zero. This happens for the unique circulant core Hadamard matrix of
order $36$~\cite{KKS:2006}.

}
\end{rem}

{\bf Complexity:} Roughly of order $v^32^k$. In this count $v^3$ is the complexity of computing an echelonized basis of $H-\sqrt{v}I$. The complexity of the invertible minor finding
algorithm is of the same order or less.
\subsection{Groebner bases}

The system $\mathcal{H}X=X$ with $X \in \{\pm 1\}^v$ can be thought of as the real quadratic system $\mathcal{H}X=X,\, \forall i \in [1,v],\, X_i^2=1$.
For background material on Groebner bases we refer the reader to \cite{AL}.

More concretely, we can consider the following steps.

\begin{enumerate}
 \item [\rm(i)] Construct the ring $P$ of polynomial functions in $v$ variables $X_i, \, i=1,\dots v.$
 \item [\rm(ii)] Construct the linear constraints $\mathcal{H}X=X.$
 \item [\rm(iii)] Construct the quadratic constraints $\forall i \in [1,v],\, X_i^2=1$
 \item [\rm(iv)] Compute a Groebner basis for the ideal $I$ of  $P$ determined by constraints (ii) and (iii).
 \item [\rm(v)] Compute the solutions as the zeros determined by $I.$
\end{enumerate}

With a tip from Delphine Boucher we produced the following program in Magma~\cite{Ma} in the case $v=4$. This program is easy to adapt for higher $v$'s.
We give it here for exposition purpose only.

{\tt F:=RationalField();} \\
//Polynomial ring defining the variables

{\tt var := 4;

 P$<w,x,y,z>$ := PolynomialRing(F,var);} \\ \ \

//The equations of the system one wants to solve over F
$$sys:=[w+x+y+z-2*w,w-x+y-z-2*x,w+x-y-z-2*y,w-x-y+z-2*z, $$ $$w^2-1,x^2-1,y^2-1,z^2-1];$$
//The ideal of the relations

{\tt
I := ideal$<P | sys>;$} \\ \ \

//Computation of a Groebner basis (for the lexicographical order if no other order is specified)\\ \ \

{\tt Groebner(I:Faugere:=true);} \\ \ \

//The set of solutions, S \\ \ \

{\tt S:=Variety(I);
S;}\\
{\bf Complexity:} As is well-known \cite{AL}, the complexity of computing Groebner bases can be doubly exponential in the number of variables, that is $v$ here.

\section{Numerics}
\label{sec:num}

The following Table~\ref{tab:upper:bound} gives an upper bound on the dimension of the eigenspace attached to the eigenvalue 1 of $\mathcal{H}$. The row $\#$ gives the number of non-Sylvester Hadamard matrices  of given order in the Magma database~\cite{Ma}.


\begin{table}[!ht]
	\centering
	\caption{Hadamard matrices with different orders in Magma database.}
	\begin{tabular}{|c|c|c|c|c|c|c|c|}
		\hline
		$v$ & 4 & 16 & 36 & 64 & 100 & 144 & 196  \\ \hline
		$\#$ & 0 & 4 & 219 & 394 & 1 & 1   & 1  \\ \hline
		$\dim \le$ & -- & 7 & 4 & 3 & 2 & 1  & 2  \\ \hline
	\end{tabular}
	\label{tab:upper:bound}
\end{table}

Given how small these upper bounds are, the method of Subsection~\ref{sec:search:lin} is very successful. By using linear algebra method, we verify that there is no self-dual bent sequence in above (non-Sylvester) Hadamard matrices. In particular, the Magma database contains only one matrix for $v\in\{100,\; 144,\; 196\}$, respectively. We thus have to construct extra matrices as explained in the Appendix.
In Table~\ref{tab:self-dual}, each column corresponds to one type of matrix from the Appendix~\footnote{For the sake of computational complexity, in Table~\ref{tab:self-dual}, we focus on Hadamard matrices with dimensions of the eigenspace attached to the eigenvalue 1 of $\mathcal{H}$ smaller than 30.}.
\begin{table}[!ht]
	\centering
	\caption{Number of self-dual bent sequences in various Hadamard matrices with dimensions of the eigenspace attached to the eigenvalue 1 of $\mathcal{H}$ smaller than 30.}
	\renewcommand{\arraystretch}{1.2}
	\resizebox{1.0\textwidth}{!}{
		\begin{tabular}{|c|c|cc|cc|ccc|c|c|}
			\hline
			$v$ & 16 &
			\multicolumn{2}{c|}{36} &
			\multicolumn{2}{c|}{64} &
			\multicolumn{3}{c|}{100} &
			144 & 196
			\\ \hline
			Types & Sylvester &
			\multicolumn{1}{c|}{\begin{tabular}[c]{@{}c@{}} Bush \\ \cite{Janko:2001},\cite{JanKhar:2002} \end{tabular}} & Paley &
			\multicolumn{1}{c|}{\begin{tabular}[c]{@{}c@{}} Regular \\ \cite{CrnPav:2001} \end{tabular}} & \begin{tabular}[c]{@{}c@{}}Regular \\ by\\ Switching\end{tabular} &
			\multicolumn{1}{c|}{Regular~\cite{CES:2018}} & \multicolumn{1}{c|}{\begin{tabular}[c]{@{}c@{}}Regular\\ \cite{Pavcevic:96}\end{tabular}} & \begin{tabular}[c]{@{}c@{}}Regular\\ Menon\\ \cite{CES:2018}\end{tabular} &
			\begin{tabular}[c]{@{}c@{}}Bush\\ \cite{Crn:2007:144-66-30},\cite{Pavcevic:96}\end{tabular}  &
			\begin{tabular}[c]{@{}c@{}}Regular\\ \cite{Crn:2007:196-91-42} \end{tabular}
			\\ \hline
			$\#\, \mbox{of}\, H$ & 1 &
			\multicolumn{1}{c|}{29} & 1 &
			\multicolumn{1}{c|}{16} & 1 &
			\multicolumn{1}{c|}{115} & \multicolumn{1}{c|}{1} & 4 &
			 4 & 4
			\\ \hline
			$\#\, \mbox{of}\, X$ & 140 &
			\multicolumn{1}{c|}{64} & 204 &
			\multicolumn{1}{c|}{\begin{tabular}[c]{@{}c@{}} 2,~4,~6,\\ 12,~620 \end{tabular}} & 2 & \multicolumn{1}{c|}{\begin{tabular}[c]{@{}c@{}} 1024,\,1056,\\1152,\,1216,
			\\2336,\,3616,\\$\cdots \cdots$
			\\5312,\,6464
			\end{tabular}} & \multicolumn{1}{c|}{12}  & 924 &
			\begin{tabular}[c]{@{}c@{}}20,\\ 924,\\ 1052\end{tabular} & \begin{tabular}[c]{@{}c@{}}6864,\\ 12870\end{tabular}
			\\ \hline
	\end{tabular}}
	\label{tab:self-dual}
\end{table}

All detailed self-dual bent sequences for the above Hadamard matrices in Table~\ref{tab:self-dual} are publicly available on \texttt{Github}: \url{https://github.com/Qomo-CHENG/Hadamard_bent}.

\section{Application}
\label{sec:application}

A recent and original application of Hadamard bent sequences, first introduced in~\cite{SCGR:2021isit}, lies in Physical Unclonable Functions (PUFs).
PUFs can be viewed as the fingerprint of a circuit, which generate unique outputs because of uncontrollable technological dispersions during the manufacturing process of silicon chips.
They are employed for many security purposes like authentications~\cite{DPGV15}, cryptographic key generations~\cite{SKAZ20}, etc. A PUF usually generates a series of random bits (by feeding customized inputs) that uniquely depends on the corresponding circuit. Therefore, one of the metrics of the performance of a PUF is the entropy of the generated random bits~\cite{RSGD16}. For instance, we expect to generate cryptographic keys as randomly as possible in practice, leading to as high entropy as possible.

It is demonstrated in~\cite{RSGD16} that $v$ inputs generated from a Hadamard matrix (e.g, $v$ row vectors) can achieve the maximal entropy of $v$ bits in PUFs. Later on, as demonstrated  in~\cite{SCGR:2021isit}, bent sequences maximize the entropy of outputs when adding one more sequence (then $v+1$ sequences in total) to the Hadamard code $(v,2v)$.
When they exist, bent sequences reach the covering radius of the Hadamard code constructed from a Hadamard matrix as in \S \ref{sec:back} (Cf. Lemma \ref{CR}). The main conjecture of~\cite{SCGR:2021isit}, checked numerically for small $v$ (e.g., $v\leq 16$), is that maximizing the entropy of outputs is equivalent to adding a new codeword at distance the covering radius of the Hadamard code. In this respect, we construct various Hadamard matrices for different $v$ up to $196$ and verify the existence of self-dual bent sequences in this paper.


\section{Conclusion}
\label{sec:con}

We have considered the self-dual bent sequences attached to Hadamard matrices from the viewpoints of generation and symmetry. Our generation method based on linear algebra works especially well
when the eigenvalue $1$ of the normalized Hadamard matrix has low geometric multiplicity. For some matrices of order $100$ this method performs well, while the
Groebner basis method cannot finish. The lack of Hadamard matrices of order $>36$ in Magma database has led us to use the switching method of~\cite{CrnSv:Switching} to generate more matrices. In general,
it would be
a worthy research project  to enrich the known databases, even in the cases where complete enumeration of equivalence classes  is unfeasible. In the same vein, refining the classification
of Hadamard matrices for $v\le 28$ from equivalence to strong equivalence would be of interest.

We note that the concept of self-dual bent sequences being not invariant by Hadamard equivalence, classification of these become infeasible, even for matrices of small order. In general, classification at order $v$ would require
to consider $(v! 2^v)^2$ matrices for each orbit representative. This makes already $147456$ for $v=4$.

\appendix

\section{Appendix on Hadamard matrices}
\label{sec:append}
In this section we indicate that how to construct Hadamard matrices of orders not sufficiently covered in Magma Hadamard database.

\subsection{Order 16}
There are five Hadamard matrices in Magma database and the first one is of the Sylvester type.

\subsection{Order 36}
 Bush-type Hadamard matrices can be found in~\cite{Janko:2001}, \cite{JanKhar:2002}. More can be generated by switching~\cite{CrnSv:Switching}.
 \subsection{Order 64}
 16 regular Hadamard matrices were obtained from the symmetric $(64,28,12)$ designs constructed in~\cite{CrnPav:2001}. One regular Hadamard matrix was obtained by switching~\cite{CrnSv:Switching}.
\subsection{Order 100}
Two Hadamard matrices can be obtained from symmetric designs $(100,45,20)$ constructed in \cite{CES:2018}. One Hadamard matrix can be obtained
from the symmetric design $(100,45,20)$ obtained in the reference \cite{Pavcevic:96}.

The switching method described in~\cite{CrnSv:Switching} applied to the Janko-Kharaghani-Tonchev symmetric $(100,45,20)$ design of~\cite{JKT:2001} corresponding to a Bush-type Hadamard matrix of order $100$ gives
$2^{10}$ designs (including the original one), $208$ of them are pairwise non-isomorphic. These $208$ symmetric $(100,45,20$) designs give rise to $120$ pairwise non-equivalent regular Hadamard matrices. In particular, $115$ of them have dimensions of the eigenspace attached to the eigenvalue $1$ of $\mathcal{H}$ smaller than $27$.
\subsection{Order 144}
Four Bush-type Hadamard matrices can be constructed from the symmetric $(144,66,30)$ designs  in~\cite{Crn:2007:144-66-30}, \cite{Pavcevic:96}. Note that one of them has dimensions of the eigenspace attached to the eigenvalue $1$ of $\mathcal{H}$ equal to $28$.

\subsection{Order 196}
Four regular Hadamard matrices were obtained from the symmetric (196,91,42) designs built in~\cite{Crn:2007:196-91-42},
and additional two regular Hadamard matrices were obtained from the symmetric (196,91,42) designs constructed in~\cite{Crn:2005:Menon}. However, the latter two have dimensions of eigensapce attached to the eigenvalue $1$ of $\mathcal{H}$ equal to $33$, so we omit them in Table~\ref{tab:self-dual}.


\begin{thebibliography}{10}
\providecommand{\url}[1]{#1}
\csname url@samestyle\endcsname
\providecommand{\newblock}{\relax}
\providecommand{\bibinfo}[2]{#2}
\providecommand{\BIBentrySTDinterwordspacing}{\spaceskip=0pt\relax}
\providecommand{\BIBentryALTinterwordstretchfactor}{4}
\providecommand{\BIBentryALTinterwordspacing}{\spaceskip=\fontdimen2\font plus
\BIBentryALTinterwordstretchfactor\fontdimen3\font minus
  \fontdimen4\font\relax}
\providecommand{\BIBforeignlanguage}[2]{{%
\expandafter\ifx\csname l@#1\endcsname\relax
\typeout{** WARNING: IEEEtranS.bst: No hyphenation pattern has been}%
\typeout{** loaded for the language `#1'. Using the pattern for}%
\typeout{** the default language instead.}%
\else
\language=\csname l@#1\endcsname
\fi
#2}}
\providecommand{\BIBdecl}{\relax}
\BIBdecl

\bibitem{AL}
Adams, William W.; Loustaunau, Philippe (1994). \emph{An Introduction to Gr\"obner Bases.} Graduate Studies in Mathematics. Vol. 3. 

\bibitem{BJH:99}
T.~Beth, D.~Jungnickel, and H.~Lenz, \emph{\BIBforeignlanguage{English}{Design
  Theory. Vol. I.}}, ser. Encycl. Math. Appl.\hskip 1em plus 0.5em minus
  0.4em\relax Cambridge: Cambridge University Press, 1999, vol.~69,
  \DOI{10.1017/CBO9780511549533}.

\bibitem{Biggs:74}
N.~Biggs, \emph{\BIBforeignlanguage{English}{Algebraic Graph Theory}}, ser.
  Camb. Tracts Math.\hskip 1em plus 0.5em minus 0.4em\relax Cambridge:
  Cambridge University Press, 1974, vol.~67.

\bibitem{Ma}
\BIBentryALTinterwordspacing
W.~Bosma, J.~J. Cannon, C.~Fieker, and A.~Steel, Eds., \emph{Handbook of Magma
  functions, Edition 2.16}, 2010. [Online]. Available:
  \url{http://magma.maths.usyd.edu.au/magma/handbook/}
\BIBentrySTDinterwordspacing

\bibitem{Bush:71}
K.~A. Bush, ``\BIBforeignlanguage{English}{Unbalanced {H}adamard matrices and
  finite projective planes of even order},''
  \emph{\BIBforeignlanguage{English}{\href{http://www.sciencedirect.com/science/journal/00973165}{J.
  Comb. Theory, Ser.~A}}}, vol.~11, no.~1, pp. 38--44, 1971,
  \DOI{10.1016/0097-3165(71)90005-7}.

\bibitem{CDPS:2010}
C.~Carlet, L.~E. Danielsen, M.~G. Parker, and P.~Sol\'e,
  ``\BIBforeignlanguage{English}{Self-dual bent functions},''
  \emph{\BIBforeignlanguage{English}{Int. J. Inf. Coding Theory}}, vol.~1,
  no.~4, pp. 384--399, Apr. 2010, \DOI{10.1504/IJICOT.2010.032864}.

\bibitem{Crn:2005:Menon}
\BIBentryALTinterwordspacing
D.~Crnkovi\'c, ``\BIBforeignlanguage{English}{Some new {M}enon designs with
  parameters $(196,91,42)$},''
  \emph{\BIBforeignlanguage{English}{\href{http://www.mathos.hr/mc/}{Math.
  Commun.}}}, vol.~10, no.~2, pp. 169--175, Dec. 2005. [Online]. Available:
  \url{https://hrcak.srce.hr/647}
\BIBentrySTDinterwordspacing

\bibitem{Crn:2006}
D.~Crnkovi\'c, ``\BIBforeignlanguage{English}{A series of regular {H}adamard
  matrices},''
  \emph{\BIBforeignlanguage{English}{\href{http://link.springer.com/journal/10623}{Des.
  Codes Cryptography}}}, vol.~39, no.~2, pp. 247--251, May 2006,
  \DOI{10.1007/s10623-005-3634-3}.

\bibitem{Crn:2007:144-66-30}
D.~Crnkovi\'c, ``\BIBforeignlanguage{English}{A construction of some symmetric
  $(144,66,30)$ designs},'' \emph{\BIBforeignlanguage{English}{J. Appl. Algebra
  Discrete Struct.}}, vol.~5, no.~1, pp. 33--39, 2007.

\bibitem{Crn:2007:196-91-42}
D.~Crnkovi\'c, ``\BIBforeignlanguage{English}{A construction of some symmetric designs
  with parameters $(196,91,42)$},'' \emph{\BIBforeignlanguage{English}{Int.
  Math. Forum}}, vol.~2, no. 61-64, pp. 3021--3026, 2007,
  \DOI{10.12988/imf.2007.07275}.

\bibitem{CES:2018}
D.~{Crnkovi\'c}, R.~{Egan}, and A.~{\v{S}vob},
  ``\BIBforeignlanguage{English}{Orbit matrices of {H}adamard matrices and
  related codes},''
  \emph{\BIBforeignlanguage{English}{\href{http://www.sciencedirect.com/science/journal/0012365X}{Discrete
  Math.}}}, vol. 341, no.~5, pp. 1199--1209, May 2018,
  \DOI{10.1016/j.disc.2018.01.018}.

\bibitem{CrnPav:2001}
D.~{Crnkovi\'c} and M.-O. {Pav\v{c}evi\'c}, ``\BIBforeignlanguage{English}{Some
  new symmetric designs with parameters (64, 28, 12)},''
  \emph{\BIBforeignlanguage{English}{\href{http://www.sciencedirect.com/science/journal/0012365X}{Discrete
  Math.}}}, vol. 237, no. 1-3, pp. 109--118, June 2001,
  \DOI{10.1016/S0012-365X(00)00364-2}.

\bibitem{CrnSv:Switching}
\BIBentryALTinterwordspacing
D.~Crnkovi\'c and A.~\v{S}vob, ``Switching for $2$-designs,''
  \emph{\BIBforeignlanguage{English}{\href{http://link.springer.com/journal/10623}{Des.
  Codes Cryptography}}},
  \DOI{10.1007/s10623-022-01059-7}. 
\BIBentrySTDinterwordspacing

\bibitem{deLauFla:2011}
W.~de~Launey and D.~Flannery, \emph{\BIBforeignlanguage{English}{Algebraic
  Design Theory}}, ser. Math. Surv. Monogr.\hskip 1em plus 0.5em minus
  0.4em\relax Providence, RI: American Mathematical Society (AMS), 2011, vol.
  175, \DOI{10.1090/surv/175}.

\bibitem{deLauSta:2008}
W.~{de Launey} and R.~M. {Stafford}, ``\BIBforeignlanguage{English}{On the
  automorphisms of {P}aley's type {II} {H}adamard matrix},''
  \emph{\BIBforeignlanguage{English}{\href{http://www.sciencedirect.com/science/journal/0012365X}{Discrete
  Math.}}}, vol. 308, no.~13, pp. 2910--2924, July 2008,
  \DOI{10.1016/j.disc.2007.07.118}.

\bibitem{DPGV15}
J.~Delvaux, R.~Peeters, D.~Gu, and I.~Verbauwhede, ``A survey on lightweight
  entity authentication with strong {PUF}s,'' \emph{ACM Comput. Surv.},
  vol.~48, no.~2, pp. Article No. 26 (1--42), Nov. 2015, \DOI{10.1145/2818186}.

\bibitem{FSSW:2013}
T.~{Feulner}, L.~{Sok}, P.~{Sol\'e}, and A.~{Wassermann},
  ``\BIBforeignlanguage{English}{Towards the classification of self-dual bent
  functions in eight variables},''
  \emph{\BIBforeignlanguage{English}{\href{http://link.springer.com/journal/10623}{Des.
  Codes Cryptography}}}, vol.~68, no. 1-3, pp. 395--406, Sept. 2013,
  \DOI{10.1007/s10623-012-9740-0}.

\bibitem{Hall:62:Mathieu}
M.~Hall, Jr, ``\BIBforeignlanguage{English}{Note on the mathieu group
  {$M_{12}$}},'' \emph{\BIBforeignlanguage{English}{Arch. Math.}}, vol.~13, pp.
  334--340, Dec. 1962, \DOI{10.1007/BF01650080}.

\bibitem{Horadam:2007}
K.~J. Horadam, \emph{\BIBforeignlanguage{English}{{H}adamard Matrices and Their
  Applications}}.\hskip 1em plus 0.5em minus 0.4em\relax Princeton, NJ:
  Princeton University Press, 2007.

\bibitem{Janko:2001}
Z.~Janko, ``\BIBforeignlanguage{English}{The existence of a {B}ush-type
  {H}adamard matrix of order $36$ and two new infinite classes of symmetric
  designs},''
  \emph{\BIBforeignlanguage{English}{\href{http://www.sciencedirect.com/science/journal/00973165}{J.
  Comb. Theory, Ser.~A}}}, vol.~95, no.~2, pp. 360--364, Aug. 2001,
  \DOI{10.1006/jcta.2000.3166}.

\bibitem{JanKhar:2002}
Z.~{Janko} and H.~{Kharaghani}, ``\BIBforeignlanguage{English}{A block
  negacyclic {B}ush-type {H}adamard matrix and two strongly regular graphs},''
  \emph{\BIBforeignlanguage{English}{\href{http://www.sciencedirect.com/science/journal/00973165}{J.
  Comb. Theory, Ser.~A}}}, vol.~98, no.~1, pp. 118--126, Apr. 2002,
  \DOI{10.1006/jcta.2001.3231}.

\bibitem{JKT:2001}
Z.~Janko, H.~Kharaghani, and V.~D. Tonchev,
  ``\BIBforeignlanguage{English}{{B}ush-type {H}adamard matrices and symmetric
  designs},''
  \emph{\BIBforeignlanguage{English}{\href{http://onlinelibrary.wiley.com/journal/10.1002/(ISSN)1520-6610}{J.
  Comb. Des.}}}, vol.~9, no.~1, pp. 72--78, Jan. 2001,
  \DOIURL{10.1002/1520-6610(2001)9:1$<$72::AID-JCD6$>$3.0.CO;2-M}{10.1002/1520-6610(2001)9:1<72::AID-JCD6>3.0.CO;2-M}.

\bibitem{Janusz:2007}
G.~J. Janusz, ``\BIBforeignlanguage{English}{Parametrization of self-dual codes
  by orthogonal matrices},''
  \emph{\BIBforeignlanguage{English}{\href{http://www.sciencedirect.com/science/journal/10715797}{Finite
  Fields Appl.}}}, vol.~13, no.~3, pp. 450--491, July 2007,
  \DOI{10.1016/j.ffa.2006.05.001}.

\bibitem{Kantor:69:aut}
W.~M. Kantor, ``\BIBforeignlanguage{English}{Automorphism groups of {H}adamard
  matrices},''
  \emph{\BIBforeignlanguage{English}{\href{http://www.sciencedirect.com/science/journal/00219800}{J.
  Comb. Theory}}}, vol.~6, no.~3, pp. 279--281, Apr. 1969,
  \DOI{10.1016/S0021-9800(69)80088-8}.

\bibitem{Kharaghani:2000:twin}
H.~Kharaghani, ``\BIBforeignlanguage{English}{On the twin designs with the
  {I}onin-type parameters},''
  \emph{\BIBforeignlanguage{English}{\href{http://www.combinatorics.org}{Electr.
  J. Comb.}}}, vol.~7, no. \#R1(1-11), 2000, \DOI{10.37236/1479}.

\bibitem{KKS:2006}
I.~S. Kotsireas, C.~Koukouvinos, and J.~Seberry,
  ``\BIBforeignlanguage{English}{{H}adamard ideals and {H}adamard matrices with
  circulant core},'' \emph{\BIBforeignlanguage{English}{J. Comb. Math. Comb.
  Comput.}}, vol.~57, pp. 47--63, 2006.

\bibitem{Kutsenko:2020}
A.~Kutsenko, ``\BIBforeignlanguage{English}{The group of automorphisms of the
  set of self-dual bent functions},''
  \emph{\BIBforeignlanguage{English}{\href{http://link.springer.com/journal/volumesAndIssues/12095}{Cryptogr.
  Commun.}}}, vol.~12, no.~5, pp. 881--898, June 2020,
  \DOI{10.1007/s12095-020-00438-y}.

\bibitem{LMS:2006}
K.~H. {Leung}, S.~L. {Ma}, and B.~{Schmidt}, ``\BIBforeignlanguage{English}{New
  {H}adamard matrices of order $4p^2$ obtained from {J}acobi sums of order
  $16$},''
  \emph{\BIBforeignlanguage{English}{\href{http://www.sciencedirect.com/science/journal/00973165}{J.
  Comb. Theory, Ser.~A}}}, vol. 113, no.~5, pp. 822--838, July 2006,
  \DOI{10.1016/j.jcta.2005.07.011}.

\bibitem{MWS}
F.~J. MacWilliams and N.~J.~A. Sloane, \emph{The Theory of Error-Correcting
  Codes}.\hskip 1em plus 0.5em minus 0.4em\relax Amsterdam, Netherlands: North
  Holland, 1977.

\bibitem{Mesnager:2016:Bent}
S.~Mesnager, \emph{\BIBforeignlanguage{English}{Bent Functions. {F}undamentals
  and Results}}.\hskip 1em plus 0.5em minus 0.4em\relax Cham: Springer, 2016,
  \DOI{10.1007/978-3-319-32595-8}.

\bibitem{Pavcevic:96}
\BIBentryALTinterwordspacing
M.-O. Pav\v{c}evi\'c, ``\BIBforeignlanguage{English}{Symmetric designs of
  {M}enon series admitting an action of {F}robenius groups},''
  \emph{\BIBforeignlanguage{English}{Glas. Mat., III. Ser.}}, vol.~31, no.~2,
  pp. 209--223, Dec. 1996. [Online]. Available:
  \url{http://books.google.com/books?id=wdgsTPYo92YC&pg=PA209}
\BIBentrySTDinterwordspacing

\bibitem{RSGD16}
O.~Rioul, P.~Sol\'e, S.~Guilley, and J.-L. Danger, ``On the entropy of
  physically unclonable functions,'' in \emph{IEEE International Symposium on
  Information Theory, Barcelona, Spain, July 10--15, 2016}.\hskip 1em plus
  0.5em minus 0.4em\relax IEEE, 2016, pp. 2928--2932,
  \DOI{10.1109/ISIT.2016.7541835}.

\bibitem{SKAZ20}
A.~Shamsoshoara, A.~Korenda, F.~Afghah, and S.~Zeadally, ``A survey on physical
  unclonable function {(PUF)}-based security solutions for {I}nternet of
  {T}hings,'' \emph{Computer Networks}, vol. 183, p. 107593, Dec. 2020,
  \DOI{10.1016/j.comnet.2020.107593}.

\bibitem{SCGR:2021isit}
P.~Sol\'e, W.~Cheng, S.~Guilley, and O.~Rioul, ``Bent sequences over {H}adamard
  codes for physically unclonable functions,'' in \emph{IEEE International
  Symposium on Information Theory, Melbourne, Australia, July 12--20,
  2021}.\hskip 1em plus 0.5em minus 0.4em\relax IEEE, 2021, pp. 801--806,
  \DOI{10.1109/ISIT45174.2021.9517752}.

\bibitem{Tokareva:2015}
N.~Tokareva, \emph{\BIBforeignlanguage{English}{Bent Functions. {R}esults and
  Applications to Cryptography}}.\hskip 1em plus 0.5em minus 0.4em\relax
  Amsterdam: Elsevier/Academic Press, 2015, \DOI{10.1016/C2014-0-02922-X}.

\bibitem{WSW:72}
W.~D. {Wallis}, A.~P. {Street}, and J.~S. {Wallis},
  \emph{\BIBforeignlanguage{English}{Combinatorics: Room Squares, Sum-Free
  Sets, {H}adamard Matrices}}, ser. Lect. Notes Math.\hskip 1em plus 0.5em
  minus 0.4em\relax Berlin: Springer-Verlag, 1972, vol. 292,
  \DOI{10.1007/BFb0069907}.

\end{thebibliography}

\providecommand\href[2]{#2} \providecommand\url[1]{\href{#1}{#1}}
  \def\DOI#1{{\small {DOI}:
  \href{http://dx.doi.org/#1}{#1}}}\def\DOIURL#1#2{{\small{DOI}:
  \href{http://dx.doi.org/#2}{#1}}}

\end{document}